\documentclass[leqno]{amsart}

\usepackage{amscd,amssymb,amsmath}

\usepackage[bookmarks]{hyperref}
\hypersetup{backref, colorlinks=true}
\usepackage{tikz}

 \usepackage{graphicx}
\usepackage[all]{xy}
\newtheorem{theorem}{Theorem}[section]

\newtheorem{proposition}[theorem]{Proposition}
\newtheorem{corollary}[theorem]{Corollary}
\newtheorem{remark}[theorem]{Remark}

\newtheorem{example}[theorem]{Example}
\newtheorem{definition}[theorem]{Definition}

\def\C{ \mathbb{C}}

\def\R{ \mathbb{R}}

\def\N{ \mathbb{N}}

\def\K{ \mathbb{K}}

\def\rond{\mathaccent"7017}

\begin{document}

\title[]{Some classes satisfying the 2-dimensional Jacobian conjecture and a proof of the complex conjecture until degree 104} 
\makeatother

\author[Nguyen Thi Bich Thuy]{Thuy Nguyen - S\~ao Paulo State Univesity (UNESP), Brazil}
\address[Nguyen Thi Bich Thuy]{S\~ao Paulo State Univesity (UNESP), S\~ao Jos\'e do Rio Preto, Brazil}
\email{bich.thuy@unesp.br}
\maketitle \thispagestyle{empty}

\begin{abstract}
We construct a non-proper set of  two variables polynomial maps and study the nowhere vanishing Jacobian condition 
 of the Jacobian conjecture for this set. 
 We obtain some classes of polynomial maps satisfying the 2-dimensional Jacobian conjecture for both real and complex cases. 
 In addition, by Newton polygon technique, we prove that the complex conjecture is true until degree 104, improving Moh boundary (degree 100) since 1983. 
 \end{abstract}

\section{Introduction} 
Let $F=(F_1,F_2): \K^2 \to \K^2$, where $\K = \R$ or $\C$, be a polynomial map satisfying the condition: 
\begin{equation} \label{Jacobian-condition}
\det (JF(x,y)) = {\rm constant} \neq 0, \quad \forall (x,y) \in \K^2,
\end{equation}
where $JF(x,y)$ is the Jacobian matrix of $F$ at $(x,y)$. 
The famous Jacobian conjecture, stated by O. H. Keller \cite{kel} in 1939, saying that 
any polynomial map $F: \K^2 \to \K^2$ with the condition  (\ref{Jacobian-condition}) 
is an automorphism, is still open today when $\K = \C$. 
In the real case, the conjecture was solved negatively  by Pinchuk \cite{Pinchuk}  in 1994. 
In this paper, following the book \cite{Essen} we call the condition (\ref{Jacobian-condition}) {\it Jacobian condition}. 

If $F$ has a global inverse, then its inverse, being continous, maps compact sets into compact sets, in other words $F$ is proper. 
The smallest set $S_F$ such that the map $F:\K^2 \, \setminus \, F^{-1}(S_F)\to \K^2 \, \setminus \, S_F $ is proper 
is called the asymptotic set of $F$. 
The Jacobian conjecture  reduces to show that the asymptotic set of  a polynomial  map satisfying the Jacobian  condition is empty. 
For polynomial maps, this set can be expressed as the set of points in the target space such that there is a sequence $\{z_k\} \subset \K^2$ tending to infinity but  its image $F(z_k)$ converges to a finite point: 
$$S_F = \{ a \in \K^2: \text{ there is a sequence }  z_k \to \infty,  F(z_k) \to a\}.$$

This paper is inspired by the counterexamples of Pinchuk \cite{Pinchuk}: 
given $(x,y)\in \mathbb{R}^2$, denote
$$t= xy -1, \quad h=t(xt+1), \quad f=(xt+1)^2(t^2+y),$$
then the Pinchuk maps $F=(p, q)$ are the ones with $p= f+h$ and $q$ varies 
 for different maps $(p,q)$ but $q$ always has the form 
$$q=-t^2 -6th(h+1) -u(f,h)$$
where $u$ is an auxiliary polynomial in $f$ and $h$ and is chosen such that 
$${\rm det} (J(p,q)) = t^2 +f^2 + (t+ f(13+15h))^2$$
(\cite[Lemma 2.2]{Pinchuk}). 
Then ${\rm det} (J(p,q)(x,y)) >0$ for every $(x,y) \in \mathbb{R}^2$ since if $t = 0$ then $f= y = \frac{1}{x} \neq 0$.  
 However, $(p, q)$ is not automorphism. 
 Our first observation is:  
 by changing variables $t:=t(x,y)$, $h:=h(x,y)$ and $f:=f(x,y)$, we can write 
 the Pinchuk maps $(p, q)$ as polynomial maps of new variables $t, h$ and $f$. 
 The sequence $\{z_k\}=\{(k,1/k)\} \subset \R^2$ tends to infinity but $t(z_k)=h(z_k)=0$ and $f(z_k)=1/k$.  
 Then the Pinchuk maps $(p, q)$, being polynomial maps of new variables $t, h$ and $f$, 
 satisfy that $(p(z_k), q(z_k))$ tends to a finite point and  therefore they are non-proper and hence, is not automorphism. 
 We call  the set of such new variables $\{ t, h, f \}$ {\it a set of non-proper variables} of the map $F=(p, q)$. 
By studying the Jacobian condition  for polynomial maps of these new variables and classifying them, 
we obtain some classes that satisfy the Jacobian conjecture for both real and complex cases.  
The idea of ``non-proper variables'' first appeared in \cite{ThuyThesis}, where we called them ``permanent variables''. We have used this idea to classify and stratify the non-proper set $S_F$ in order to calculate the intersection homology of singular varieties constructed  in \cite{ThuyValette} (see also in \cite{ThuyCidinha} and \cite{Thuy_pinchuk}). This is another our approach to the Jacobian conjecture, by intersection homology. In this article, we also use Newton polygon technique in \cite{Abhyankar, Nagata2, Nagata3, Appelgate-Onishi, Magnus, Nakia-Baba, Henryk} 
to prove that the 2-dimensional complex Jacobian conjecture is true until degree 104, improving Moh boundary \cite{Moh} (1983) from 100 till 104. 

We would like to remark that in this article we try to provide 
another point of view to investigate the Jacobian conjecture using algebraical geometrical and topological approach, 
but there have been many approaches to the conjecture. 
One of the newest approaches is via quantization \cite{Alex, BK}, relating the Jacobian conjecture with the Dixmier conjecture.

 \section*{Structure and Results of the paper}
 
 The paper has two parts. 
 
 The first part, presented in Section  \ref{pertinent-variables}, 
 concerns the construction of non-proper variables for both real and complex cases and a classification of non-proper maps under the Jacobian condition. We obtain some  explicit classes satisfying the 2-dimensional Jacobian conjecture for both real and complex cases (Corollaries \ref{cor:2 p-r=1}, \ref{cor n=1 Comp} and \ref{cor n=1 Real}). These results are classified in Remark \ref{remark classes}. 
 We also provide a possible counter-example for the conjecture (if there exists) in Corollary \ref{cor-counterexample}. 

 The second part, presented in  Section \ref{degree105-complex}, concerns a particular result to the complex  conjecture for degree until 104 (Theorem  \ref{theorem-degree104}). We prove that the complex conjecture is true for degree until 104 
using known results studied by 
Abhyankar \cite{Abhyankar}, Nagata \cite{Nagata2, Nagata3}, 
 Appelgate and Onishi  \cite{Appelgate-Onishi}, Magnus \cite{Magnus}, 
Nakai and Baba \cite{Nakia-Baba} and 
\.{Z}o\l\c{a}dek \cite{Henryk}. 
It is known that the complex conjecture is true for degree until 100, proved by Moh \cite{Moh} in 1983 with a highly non-trivial proof using computer assistant improvement. 
 There are also some works increasing the degree 100 of Moh, for instance \cite {GGHV}, but we would like to emphasize 
that our proof is not a very difficult one, that may be surprising for this mysterious conjecture.  

Let us remark that there may be a connection between the two parts of the paper to study the complex conjecture by this direction for higher degree cases (Remark \ref{Remark-connection}).

\medskip

\section{Non-proper variables sets of a polynomial map} \label{pertinent-variables}

\subsection{Construction and examples}

We consider $\K = \R$ or $\C$ and $F=(F_1, F_2): \K^2 \to \K^2$  a polynomial map of two variables $x$ and $y$. 
 With $z \in \K^2$, we denote  by $\|z\|$ its Euclidean norm. 
Let us recall firstly that the polynomial map $F$ is non-proper if there is a sequence 
$\{z_k = (x_k, y_k)\}_{k\in\N}\subset \K^2$ such that $\lim_{k\to\infty} \|z_k\|= + \infty$ and  $\lim_{k\to\infty}F(z_k)=(c_1, c_2)$, for some finite $c_1, c_2 \in \K$. 
The following remark is important for the definition of  {\it non-proper variables sets} 
that we will introduce in Definition \ref{d:p-v}. 
\begin{remark}
 $F$ is non-proper if and only if it can be written as a polynomial map of some variables $u_0, u_1, \dots, u_n$ ($\, n \geq 1$), where $u_0, u_1, \dots, u_n$ are polynomial functions of two variables $x$ and $y$ such that there is a sequence 
$\{z_k= (x_k, y_k)\} \subset \K^2$ with 
$\lim_{k\to\infty} \|z_k\|= + \infty$ but 
 $\lim_{k\to\infty} u_i(z_k)=\alpha_i \in \K$ for all $\alpha_i$ are finite numbers, i.e.
\[F(x,y)=H(u_0(x,y), u_1(x, y), \ldots, u_n(x,y)), \]
where $H=(H_1, H_2): \K^{n+1} \to \K^2$ is a polynomial map of variables $u_0, u_1, \dots, u_n$. 
 \end{remark}
 \noindent In fact, if $F = (F_1, F_2)$ is non-proper, we can take $u_0:=F_1$ and $u_1:=F_2$ and the inverse is obvious. 
 We define:

\begin{definition}[non-proper variables sets] \label{d:p-v} 
{\rm 
Let $F=(F_1, F_2): \K^2 \to \K^2$  be a polynomial map. 
A set of polynomial functions $u_0, u_1, \dots, u_n: \K^2 \to \K \, (n \geq 1)$  
is called a set of 
{\it non-proper variables}  of $F$ if there exists 
\begin{enumerate} 
\item a sequence $\{z_k \}_{k\in\N}\subset \K^2$ such that 
$\lim_{k\to\infty} \|z_k\|= + \infty$ and  $\lim_{k\to\infty} u_i(x_k,y_k)$ is a finite number $\alpha_i \in \K$ for all $i= 0, \dots, n$;
\item a polynomial map $H=(H_1, H_2): \K^{n+1} \to \K^2$ such that 
$$F(x,y)=(H_1(u_0, u_1, \ldots, u_n), H_2(u_0, u_1, \ldots, u_n)).$$
\end{enumerate} 
}
\end{definition} 

\begin{remark} \label{remark-non-proper-map}
{\rm 
A polynomial map $F: \K^2 \to \K^2$  is non-proper if and only if it possesses at least a set of non-proper variables $\{u_0, u_1, \dots, u_n\}$.
}
\end{remark}

\begin{definition}[Dependent/Independent non-proper variables sets] \label{d:dep/inde}
{\rm 	 

Let $F: \K^2 \to \K^2$ be a non-proper polynomial map and $\{u_0, u_1, \dots, u_n\}$  is a set of non-proper variables of $F$.  
 We say that: 
 \begin{enumerate} 
 \item $\{u_0, u_1, \dots, u_n\}$ is 
{\it dependent} if there exists $z=(x,y) \in \K^2$ such that the two vectors 
	$$\left( \frac{\partial u_0}{\partial x} (z), \frac{\partial u_1}{\partial x}(z), \ldots, \frac{\partial u_n}{\partial x}(z) \right), $$
	$$\left(\frac{\partial u_0}{\partial y} (z), \frac{\partial u_1}{\partial y}(z), \ldots, \frac{\partial u_n}{\partial y}(z) \right)$$
	are linearly dependent in $\K^{n+1}$; 
\item  otherwise, $\{u_0, u_1, \dots, u_n\}$ is {\it independent} if for any $z=(x,y)\in\K^2$, the two vectors 
	$$\left( \frac{\partial u_0}{\partial x} (z), \frac{\partial u_1}{\partial x}(z), \ldots, \frac{\partial u_n}{\partial x}(z) \right), $$
	$$\left( \frac{\partial u_0}{\partial y}(z), \frac{\partial u_1}{\partial y}(z), \ldots, \frac{\partial u_n}{\partial y}(z) \right)$$
	are linearly independent in $\K^{n+1}$.
\end{enumerate}	
	}	
\end{definition}  

\begin{example}[{\cite[Section 6]{Chau}}] 
{\rm Let $F: \C^2 \to \C^2$ be the polynomial map defined by: 
$$F(x,y) = (x - 2(xy+1) - y(xy+1)^2, -1 - y(xy +1)).$$ 
Put: 
$$u_0 = x, \quad u_1= xy +1, \quad u_2 = y(xy+1), \quad u_3 = y(xy+1)^2, $$ 
then $\{u_0, u_1, u_2, u_3\}$ is a set of non-proper variables of $F$ 
since any sequence $\{ (1/k, -k+ c)\}$, 
with  $c \in \C$, satisfies the first condition in Definition \ref{d:p-v}. 
The second condition is verified by taking 
$$H: \C^4 \to \C^2, \qquad H(u_0, u_1, u_2, u_3) = (u_0 - 2u_1 -u_3, -1 - u_2)$$ 
and then $F = H(u_0, u_1, u_2, u_3)$. 
Consequently, $F$ is non-proper.

Now we show that $\{u_0, u_1, u_2, u_3\}$ is independent. In fact, for any $z=(x,y) \in \C^2$, we have
\[
\begin{matrix}
\left( \frac{\partial u_0}{\partial x} (z), \frac{\partial u_1}{\partial x}(z), \frac{\partial u_2}{\partial x}(z), \frac{\partial u_3}{\partial x}(z) \right)  
& = & (1 , &y, &y^2, &2xy^3 + 2y^2), \cr 
\left( \frac{\partial u_0}{\partial y}(z), \frac{\partial u_1}{\partial y}(z), \frac{\partial u_2}{\partial y}(z), \frac{\partial u_3}{\partial y}(z) \right) 
& = & (0, & x, & 2xy +1, &3x^2y^2 + 4xy+1). 
\end{matrix}
\]
We see that there is no point $z \in \C^2$ such that the second partial derivative at $z$ is the zero vector in $\C^4$. Then the two partial derivatives above are linearly independent vectors at any point $z \in \C^2$. 
By Definition \ref{d:dep/inde}, the set of non-proper variables $\{u_0, u_1, u_2, u_3\}$ is independent.

 Let us remark that if we consider the above map as a real one $F: \R^2 \to \R^2$, i.e. the same map $F$ but with real variables, then $\{u_0, u_1, u_2, u_3\}$ is also an independent set of non-proper variables.}
\end{example}

\begin{proposition} \label{p:2}  Let $F=(F_1, F_2):\K^2 \to\K^2$ be a non-proper polynomial map and 
$\{u_0, u_1, \dots, u_n\}$  be a set of non-proper variables of $F$.  
If $\{u_0, u_1, \ldots, u_n\}$ is dependent  
 then $F$ does not satisfy the Jacobian condition.  
\end{proposition} 

\begin{proof} 
By Definition  \ref{d:p-v}, there exists a polynomial map 
$H=(H_1,H_2):\K^{n+1}\to\K^2$ such that 
	$$F_1(x,y) = H_1(u_0(x,y), u_1(x,y), \ldots, u_n(x,y)),$$ 
	$$F_2(x,y) = H_2(u_0(x,y), u_1(x,y), \ldots, u_n(x,y)).$$
We have 
$$
    JF(x,y) =     \begin{pmatrix}
            \frac{\partial F_1}{\partial x} & \frac{\partial F_1}{\partial y} \cr 
            \frac{\partial F_2}{\partial x} & \frac{\partial F_2}{\partial y}
        \end{pmatrix} 
    = \begin{pmatrix}
           \frac{\partial H_1}{\partial u_0} & \frac{\partial H_1}{\partial u_1} & \cdots & \frac{\partial H_1}{\partial u_n} \cr 
           \frac{\partial H_2}{\partial u_0} & \frac{\partial H_2}{\partial u_1} & \cdots & \frac{\partial H_2}{\partial u_n}
           \end{pmatrix} 
      \begin{pmatrix}
          \frac{\partial u_0}{\partial x} & \frac{\partial u_0}{\partial y} \cr 
          \frac{\partial u_1}{\partial x} & \frac{\partial u_1}{\partial y} \cr 
          \vdots & \vdots \cr 
          \frac{\partial u_n}{\partial x} & \frac{\partial u_n}{\partial y} 
       \end{pmatrix}. 
$$
Since $u_0, u_1, \ldots, u_n$ are dependent, then by Definition \ref{d:dep/inde}, 
there exists $z = (x,y) \in \K^2$ such that the  two vectors  
$$\left( \frac{\partial u_0}{\partial x}(z), \frac{\partial u_1}{\partial x}(z), \ldots, \frac{\partial u_n}{\partial x}(z) \right), $$
$$\left( \frac{\partial u_0}{\partial y}(z), \frac{\partial u_1}{\partial y}(z), \ldots, \frac{\partial u_n}{\partial y}(z) \right)$$
are linearly dependent. It follows that ${\rm det}( JF(z)) =0$. 
\end{proof}

\begin{corollary}  \label{class Real-Complex}
	 For any set of dependent non-proper variables $u_0, u_1, \ldots, u_n: \K^2 \to \K$  and for any polynomial map $H=(H_1, H_2):\K^{n+1} \to\K^2$, then $F:=H(u_0, u_1, \ldots, u_n)$ is a polynomial map $\K^2 \to \K^2$ satisfying the Jacobian conjecture. 
\end{corollary}

\begin{corollary}  \label{Nonzero&Independence}
Let $F=(F_1, F_2):\K^2 \to\K^2$ be a non-proper polynomial map and 
$\{u_0, u_1, \dots, u_n\}$  be a set of non-proper variables of $F$. 
If  $F$ satisfies the Jacobian condition then $\{u_0, u_1, \ldots, u_n\}$ is independent. 
\end{corollary} 

\medskip 

In the following example we will show that the functions $t, h, f$ defined in \cite{Pinchuk} and \cite{Filipe} 
form a set of non-proper variables  for the Pinchuk maps and it is independent. 
Let us remark that we will consider both of real and complex cases, {i.e.} we consider the 
polynomials $t, h, f:\K^2 \to\K$, with $\K=\R$ or $\C$. In the complex case, it shows that the inverse of Corollary  \ref{Nonzero&Independence} is not true.

\begin{example} \label{ex-Pinchuk}
{\rm 

Let us recall firstly the Pinchuk map defined in \cite{Pinchuk}: 
given $(x,y)\in \mathbb{R}^2$, denote
$$t= xy -1, \quad h=t(xt+1), \quad f=(xt+1)^2(t^2+y),$$
then the Pinchuk maps $(p, q)$ are the ones with $p= f+h$ and $q$ varies 
 for different maps $(p,q)$ but $q$ always has the form 
$$q=-t^2 -6th(h+1) -u(f,h)$$
where $u$ is an auxiliary polynomial in $f$ and $h$ is chosen such that 
$${\rm det} (J(p,q)) = t^2 +f^2 + (t+ f(13+15h))^2$$
(\cite[Lemma 2.2]{Pinchuk}).

As remarked in \cite{Pinchuk}, we have $\dfrac{h(x,y)+1}{x}=t^2(x,y)+y$. Thus, we may write: 
$$ f(x,y)=(xt(x,y)+1)^2(t^2(x,y)+y).$$

Considering the 
polynomials $t, h, f$ as both real and complex polynomials and the Pinchuk maps $(p,q)$ as both real and comlex polynomial mappings, {i.e.} we consider $t, h, f:\K^2 \to\K$ and $(p, q): \K^2 \to \K^2$, with $\K=\R$ or $\C$, we have the following affirmations: 

\begin{enumerate}

\item  $\{t, h, f\}$ is a set of non-proper variables in both cases $\K=\R$ and $\C$. 

 In fact, let $z_k:=(k,1/k)$, then $t(z_k)=h(z_k)=0$ and $f(z_k)=1/k$. 

\item $\{t, h, f\}$ is a set of independent non-proper variables (Definition \ref{d:dep/inde}) in both cases $\K=\R$ and $\C$. 

In fact, from \cite[Lemma 2.1]{Pinchuk}, it is enough to show that there is no $(x,y) \in\K^2$ such that 
$$h(x,y)-t(x,y)=-f(x,y)=-f+3h^2+2h=0.$$ 
Thus, the proof follows by a direct computation. 

\item There exists a polynomial map 
$H:\K^3\to\K^2$ such that $(p,q)(x,y) = H(t, h, f).$

In fact, as in \cite[page 3]{Pinchuk},  set $p = f+h$, then by the Lemma 2.2  in \cite{Pinchuk}, there exists a polynomial $u(f, h)$ such that 
$$\det  J(p,  -t^2 -6th(h+1)-u(f,h)) = t^2 + (t+ (f(13+15h))^2 +f^2.$$
 
We consider $H=(p, q):\K^3 \to\K^2$ the polynomial map defined by 
\[ p(t,h,f)=f+h,  \quad  q(t,h,f)= -t^2 -6th(h+1)-u(f,h).\]

The polynomials $p, q$ yield the polynomial map $H:\K^3\to\K^2$, where $H(t,h,f)=(p(t,h,f),q(t,h,f))$. 
 Then the Pinchuk maps are defined in  \cite{Pinchuk} by: 
$$(p,q):\K^2 \to\K^2, \qquad (p,q)(x,y) := H(t(x,y), h(x,y), f(x,y)).$$ 
\end{enumerate}
}
\end{example}

\begin{remark} \label{remark-Pinchuk-deg25}

{\rm 

As the remark at the end of the paper  \cite{Pinchuk}, 
the Pinchuk map constructed in the proof of  \cite[Lemma 2.2]{Pinchuk} has degree 40, 
where ${\rm deg}(p) = 10$ and  ${\rm deg}(q) = 40$, 
but one can reduce ${\rm deg} (q)$ to 35 by a suitable choice of the 
auxiliary polynomial $u(f,h)$. 
In fact, the polynomial $q$ varies for different Pinchuk map, 
but always has the form $q = -t^2 -6th(h+1) - u(f,h)$. 
In \cite{Essen} (page 241), van den Essen  choose
$$u(f, h) = 170fh + 91h^2 + 195fh^2 + 69h^3 + 75fh^3 + \dfrac{75}{4}h^4,$$
and in this case the degree of Pinchuk map is 25. 
This is also the one that Andrew Campbell studies in the series of his papers 
\cite{Andrew1, Andrew2, Andrew3}. 
Of course, this Pinchuk map also satisfies all the affirmations in the Example \ref{ex-Pinchuk}. 
}
\end{remark}

\begin{example}  \label{Example-Filipe}

{\rm In \cite{Francisco-Filipe, Filipe}, basing on Pinchuk's examples \cite{Pinchuk}, the authors provides new counterexamples to the real conjecture but with the possibility to reduce the degree of $q$ to 15: In \cite[Corollary 2.6]{Filipe}, they prove that if taking $p = f + h$ and 
$$q = -\dfrac{(h^3+h^2)^2}{f^2} - \dfrac{2h^2(h^3+h^2)}{f}+4h^3 + \dfrac{3h^2}{2}-5hf,$$ 
where $h$ and $f$ are the same polynomials defined in Pinchuk's examples \cite{Pinchuk}, then $(p, q)$ is also a non-injective polynomial local diffeomorphism. In this case ${\rm deg} \, p = 10$ and ${\rm deg } \, q = 15$. 

With the same polynomial $t$ defined in Pinchuk maps \cite{Pinchuk}, one has 
$$t = h - \dfrac{h^2(h+1)}{f}$$
(see the proof of \cite[Proposition 2.3]{Filipe}). It follows: 
$$q = -(h-t)^2 - 2h^2(h-t) +4h^3 + \dfrac{3h^2}{2}-5hf.$$
Then, in this case $q$ is also a polynomial of non-proper variables $t, h, f$. 
In other words, the examples of ``Pinchuk type'' constructed in \cite{Filipe} are also polynomial maps of non-proper variables $t, h$ and $f$ satisfying all the affirmations in Example \ref{ex-Pinchuk}. 
}
\end{example}

\subsection{A set of non-proper variables} \label{class-pertinent}

In this section, we consider the following set of polynomial functions $u_0, u_1, \dots, u_n: \K^2 \to \K \; (n \geq 1)$:  
\begin{equation} \label{eq:p-r}
 	u_0(x,y) =  y, \quad u_i(x,y) =  x^i - b^i x^{ir} y^{is}, \quad b \in \K\setminus \{0\}, \quad r, s \geq 1 \quad (i = 1, \dots, n).
\end{equation}

\noindent The motivation for studying this set of variables is that any polynomial written under the variables $u_0, u_1, \ldots, u_n$: 
$$\alpha u_0 + \beta u_1 + \text{ terms of higher degree} \qquad (\alpha \neq 0, \beta \neq 0)$$
  has Newton polygon containing two points (1, 0) and (0, 1) 
 (see Definition of Newton polygon in \cite{Abhyankar}, \cite{Nagata2, Nagata3}, also in \cite[page 93]{Essen}). 
In the following we will prove first that for any triple of non-zero natural numbers $(n, r, s)$ 
and for any polynomial mapping $F:\K^2 \to \K^2$ such that 
$$F(x,y)=(H_1(u_0, u_1, \ldots, u_n), H_2(u_0, u_1, \ldots, u_n)),$$ 
where $H=(H_1, H_2): \K^{n+1} \to \K^2$, 
then $\{u_0, u_1, \dots, u_n\}$ is a set of non-proper variables of $F$ (Proposition \ref{prop:p-r} for complex case and Proposition \ref{prop:p-r-Real} for real case). 
 Then we investigate the independence (Definition \ref{d:dep/inde}) for this set: 
 we consider all  the possibilities of triples $(r, s, n)$, we divide them into three cases: $r=1$ (section \ref{subsection p:r=1});  
$(n=1, \, r \geq 2)$ (section \ref{subsection p:r>=1}); 
and $(n \geq 2, r\geq 2)$ (section \ref{subsection p:rneq1}). 
We obtain some classes satisfying the Jacobian conjecture for polynomial mappings $F: \K^2 \to \K^2$, 
 both real and complex cases $\K = \R$ or $\C$ (Corollary \ref{cor n=1 Comp} and Corollary \ref{cor n=1 Real}).

\subsubsection{The set of variables (\ref{eq:p-r})  is a set of non-proper variables}  \label{class-non-proper}

\begin{proposition}\label{prop:p-r} 
For any fixed $n \geq 1$, the set of variables $u_0, u_1, \ldots, u_n: \C^2 \to \C$ in (\ref{eq:p-r}) is a set of non-proper variables  of any polynomial map 
$$F(x,y)=H(u_0, u_, \ldots, u_n): \C^2 \to \C^2$$ 
where $H=(H_1, H_2): \C^{n+1} \to \C^2$ is a polynomial map of variables $u_0, u_1, \dots, u_n$. 
\end{proposition}
\begin{proof}
 With a fixed $n \geq 1$ and given $r, s \geq 1$, let us consider the sequence $\{ (x_k, y_k)\} \subset \C^2$, with 
$$x_k = k, \qquad 
y_k = \sqrt[s]{\frac{1}{b} \frac{1}{k^{r-1}} \left( 1 + \frac{1}{k^n} \right)} \qquad (b \neq 0).$$ 
We will prove that $u_i(x_k, y_k)$ does not tend to infinity for any $i = 0, 1, \ldots, n$.  

At first, let us remark that $u_0(x_k,y_k) = y_k$ tends to $\sqrt[s]{1/b}$ if $r = 1$ and tends to 0 if $r >1$. 
Now for $i = 1, \ldots, n$, we have 
$$u_i(x,y) = x^i[1 - b^ix^{(r-1)i}y^{is}],$$ 
then, by a direct computation, it follows: 
$$u_i(x_k, y_k) = - \sum_{l=1}^i 
\begin{pmatrix} i \cr l \end{pmatrix} 
\frac{1}{k^{nl-i}},$$
where
$$\begin{pmatrix} i \cr l \end{pmatrix}  = \frac{i!}{l!(i-l)!}.$$
For $i = 1, \dots, n-1$, the sequence $u_i(x_k,y_k)$ tends to 0 since in this case $nl  - i >0$ for any  $l \geq 1$. 
Now we have 
$$u_n(x_k, y_k) = -\sum_{l = 1}^n
\begin{pmatrix} n \cr l \end{pmatrix} 
\dfrac{1}{k^{nl -n}} = 
- \begin{pmatrix} n \cr 1 \end{pmatrix}  
- \sum_{l = 2}^n \begin{pmatrix} n \cr l \end{pmatrix}  \dfrac{1}{k^{nl -n}},$$ 
 then $u_n(x_k, y_k)$ tends to  $(-n)$.  
\end{proof}

In the real case, if $s$ is even then the sequence $\{ (x_k, y_k)\}$ chosen in the above proof can be used only if $b >0$. 
Therefore we have:  

\begin{proposition} \label{prop:p-r-Real} 
If $s$ is odd, then the set of  variables $u_0, u_1, \ldots, u_n: \R^2 \to \R$ in (\ref{eq:p-r}) is a set of non-proper variables  of any polynomial map 
$F(x,y)=H(u_0, u_, \ldots, u_n): \R^2 \to \R^2$ 
where $H=(H_1, H_2): \R^{n+1} \to \R^2$ is a polynomial map of variables $u_0, u_1, \dots, u_n$.  
If $s$ is even, then the result holds when  $b >0$. 
\end{proposition}

\begin{remark}
{\rm 
Proposition  \ref{prop:p-r} and Proposition \ref{prop:p-r-Real} are still true if we change the sign of all  the terms in the non-proper variables (\ref{eq:p-r}), { i.e.} they are true for the set of variables: 
$$u_0(x,y) =  -y, \quad u_i(x,y) =  -x^i + b x^{ir} y^{is}, \quad r, s \geq 1 \quad (i = 1, \dots, n).$$
}
\end{remark}

\begin{remark}
{\rm Proposition  \ref{prop:p-r} and Proposition \ref{prop:p-r-Real} are still true if we exchange the roles of $x$ and $y$.} 
\end{remark}

Considering the variables $u_0, u_1, \dots, u_n$ in (\ref{eq:p-r}), from Remark  \ref{remark-non-proper-map}, we have the following corollaries: 

\begin{corollary} 
If $F: \C^2 \to \C^2$ can be written as a polynomial map of variables $u_0, u_1, \dots, u_n$, i.e. 
$F(x,y)=H(u_0, u_1, \ldots, u_n)$, 
where $H=(H_1, H_2): \C^{n+1} \to \C^2$ is a polynomial map of variables $u_0, u_1, \dots, u_n$, then  $F$ is non-proper. 
\end{corollary}

\begin{corollary}  
If $F: \R^2 \to \R^2$ can be written as a polynomial map of variables $u_0, u_1, \dots, u_n$, i.e. 
$F(x,y)=H(u_0, u_1, \ldots, u_n)$, 
where $H=(H_1, H_2): \R^{n+1} \to \R^2$ is a polynomial map of variables 
$u_0, u_1, \dots, u_n$ in (\ref{eq:p-r}) then  $F$ is non-proper if $s$ is odd. 
If $s$ is even, then the result holds when  $b >0$. 
\end{corollary}

\subsubsection{The set of non-proper variables  (\ref{eq:p-r}) with $r=1$} \label{subsection p:r=1}

We consider the non-proper variables $u_0, u_1, \dots, u_n: \K^2 \to \K$  in (\ref{eq:p-r}) with $r=1$, 
{i.e.} the set 
\begin{equation} \label{eq:p-r=1}
 	u_0 =  y, \quad u_i =  x^i - b^ix^{i} y^{is}, \quad s \geq 1 \quad (i = 1, \dots, n). 
\end{equation} 
Here we consider $b \neq 0$ if $\K = \C$ (cf. Proposition \ref{prop:p-r}) and $b >0$ if $\K = \R$ (cf. Proposition \ref{prop:p-r-Real}).

\begin{theorem} \label{prop1 p:r=1} 
For any fixed $n \geq 1$, the set $\{u_0, u_1, \ldots, u_n \}$ in (\ref{eq:p-r=1}) is a set of dependent  non-proper variables. 		
\end{theorem}
\begin{proof}
By Definition \ref{d:dep/inde}, we show that there exists $z \in \K^2$ such that 
$\dfrac{\partial u_i}{\partial x}(z) = 0$, for all $i = 0, 1, \dots, n$. 
In fact, we have 
	$$\frac{\partial u_0}{\partial x} = 0, \quad
	\frac{\partial u_i}{\partial x}= ix^{i-1}[1 - (b y^s)^i].$$
Then any point $z=(x_0,\sqrt[s]{1/b})$, where $x_0 \in \K$, is a solution of the system $\dfrac{\partial u_i}{\partial x} (z)= 0$, $i= 0, 1, \ldots, n$ 
(notice that we consider $b \neq 0$ if $\K = \C$ and $b >0$ if $\K = \R$). 
\end{proof}

It follows from Proposition \ref{p:2} and Theorem \ref{prop1 p:r=1} that: 

\begin{corollary} \label{cor:2 p-r=1}
Any polynomial map 
$F(x,y)=H(u_0, u_, \ldots, u_n): \K^2 \to \K^2$, 
where $H=(H_1, H_2): \K^{n+1} \to \K^2$ is a polynomial map of variables $u_0, u_1, \dots, u_n$  
in (\ref{eq:p-r=1}) satisfies the Jacobian conjecture. 
\end{corollary}

\subsubsection{The set of non-proper variables  (\ref{eq:p-r}) with $n=1$ and $r \geq 2$} \label{subsection p:r>=1} 
 Consider the case $n=1$, {i.e.} the case where the set (\ref{eq:p-r}) has only two variables: 
$$u_0 = y, \quad u_1 =  x - b x^{r} y^{s}.$$

\begin{theorem}  \label{theo n=1} 
We have: 
\begin{enumerate}
\item In the complex case $u_0, u_1: \C^2 \to \C$, the set of non-proper variables $\{u_0, u_1\}$ is dependent with any $r \geq 2$, $s \geq 1$ and 
$b \neq 0$; 
\item In the real case $u_0, u_1: \R^2 \to \R$, the set of non-proper variables $\{u_0, u_1\}$  is dependent with any $r$ odd, $s$ even and $b >0$.
\end{enumerate}

\end{theorem}

\begin{proof} 
We have ${\partial u_1}/{\partial x} = 1 - rbx^{r-1}y^s$. 
With the same argument as in the proof of Theorem \ref{prop1 p:r=1}, 
 the equation ${\partial u_1}/{\partial x} = 0$ always has a solution in the complex case. 
 In the real case, this equation has solution when $r$ is odd, $s$ is even and $b>0$. 
 \end{proof}

From Proposition \ref{p:2}, Corollary \ref{cor:2 p-r=1}, Theorem  \ref{theo n=1} and Corollary \ref{cor:2 p-r=1}, we have two  explicit classes in which the two-dimensional Jacobian conjecture is true, for complex and real cases, respectively: 
 
\begin{corollary}  \label{cor n=1 Comp}
Let $F = (F_1, F_2): \C^2 \to \C^2$ be a polynomial map. 
If $F_1$ and $F_2$ can be written as a finite sum: 
$$ \sum_{\alpha, \beta} a_{\alpha \beta}u_0^\alpha u_1^\beta, \qquad a_{\alpha \beta} \in \C,$$
where 
$$u_0 = y, \quad u_1 =  x - b x^{r} y^{s} \quad (b \neq 0),$$
then $F$ satisfies the Jacobian conjecture. 
\end{corollary}

\begin{corollary}  \label{cor n=1 Real}
Let $F = (F_1, F_2): \R^2 \to \R^2$ be a polynomial map. 
If $F_1$ and $F_2$ can be written as a finite sum: 
$$ \sum_{\alpha, \beta} a_{\alpha \beta}u_0^\alpha u_1^\beta, \qquad a_{\alpha \beta} \in \R,$$
where 
$$u_0 = y, \quad u_1 =  x - bx^{r} y^{s} \quad (r  \, {\rm odd}, \, s \, {\rm even}, \, b>0),$$
then $F$ satisfies the Jacobian conjecture. 
\end{corollary}

\subsubsection{The set of non-proper variables  (\ref{eq:p-r}) with $n\geq 2$ and $r \geq 2$}  \label{subsection p:rneq1}
 
We consider the set of non-proper variables $u_0, u_1, \dots, u_n$ in (\ref{eq:p-r})  with $r \geq 2$, {i.e.} the set 
\begin{equation} \label{eq:p-rneq1}
u_0 =   y, \quad u_i =    x^i - b^ix^{ir} y^{is}, \quad r \geq 2, \quad s \geq 1 \quad  (i = 1, \dots, n). 
\end{equation} 
Here we consider $b \neq 0$ if $\K = \C$ (cf. Propostion \ref{prop:p-r})
 and $b >0$ if $\K = \R$ (cf. Proposition \ref{prop:p-r-Real}).

\begin{proposition} \label{prop rgeq2 ngeq2}
For a fixed natural number $n \geq 2$, the set of variables $u_0, u_1, \ldots, u_n $ in \eqref{eq:p-rneq1} is a set of independent non-proper variables. 	
\end{proposition}
\begin{proof} 
For any $z = (x,y) \in \K^2$, we have: 
	$$\left( \frac{\partial u_0}{\partial x} (z), \frac{\partial u_1}{\partial x}(z), \frac{\partial u_2}{\partial x}(z), \ldots, \frac{\partial u_n}{\partial x}(z) \right) = \left( 0 , \frac{\partial u_1}{\partial x}(z), \frac{\partial u_2}{\partial x}(z),  \ldots, \frac{\partial u_n}{\partial x}(z) \right), $$
	$$\left( \frac{\partial u_0}{\partial y}(z), \frac{\partial u_1}{\partial y}(z), \frac{\partial u_2}{\partial y}(z), \ldots, \frac{\partial u_n}{\partial y}(z) \right)=\left( 1, \frac{\partial u_1}{\partial y}(z), \frac{\partial u_2}{\partial y}(z), \ldots, \frac{\partial u_n}{\partial y}(z) \right).$$ 
By  Definition \ref{d:dep/inde}, we need to show that the two above vectors 
are linearly independent for any $z\in\K^2$. 
 It is enough to show that  there is no $z\in\K^2$ such that $\dfrac{\partial u_1}{\partial x}(z) = \dfrac{\partial u_2}{\partial x}(z) = 0$. We have
	$$\frac{\partial u_1}{\partial x} = 1 - rbx^{r-1}y^s, \quad \frac{\partial u_2}{\partial x} = 2x[1 - (b\sqrt{r} x^{r-1}y^s)^2].$$
Let us suppose that  $\dfrac{\partial u_1}{\partial x}(z)= \dfrac{\partial u_2}{\partial x}(z) =0$. The  first equation implies $rbx^{r-1}y^s=1$, then $y \neq 0$. Since $r >1$, 
we have also $x \neq 0$. 
Therefore, the second equation implies $b \sqrt{r} x^{r-1}y^s=1$. 
From the two last equations, we have $rbx^{r-1}y^s= \sqrt{r}bx^{r-1}y^s$, that implies $r =1$, 
which is a contradiction because $r>1$. 
\end{proof}

\begin{remark} \label{remark classes}

{\rm 

Let us denote by $\mathcal{C}$ the set 
 of the non-proper polynomial maps $F = H(u_0, u_1, \dots, u_n): \K^2 \to \K^2$ 
 with $u_0, u_1, \dots, u_n$ defined in (\ref{eq:p-r}) 
and by  $\mathcal{C}^n_{r,s}$ the corresponding subset of $\mathcal{C}$ 
     for each fixed triple of natural numbers $(r, s, n)$.   
     The theorems \ref{prop1 p:r=1} and \ref{theo n=1} 
    say that the classes 
    $\mathcal{C}^n_{1,s}$ ($n \geq 1$, $s \geq 1$) and 
    $\mathcal{C}^1_{r,s}$ ($r \geq 1$, $s \geq 1$) satisfy the complex Jacobian conjecture. 
    This fact holds for the real conjecture if $r$ odd, $s$ even and $b >0$. 
    
    \medskip 
     
 Moreover, from Proposition  \ref{prop rgeq2 ngeq2} and Remark \ref{Nonzero&Independence}, we have: 

\begin{corollary} \label{cor-counterexample}

A polynomial map belonging to a class {\rm $\mathcal{C}^n_{r, s}$ ($n \geq 2, r\geq 2, s \geq 1$)} 
satisfying the Jacobian condition is a counterexample to the Jacobian conjecture. 
\end{corollary} 
}
    \end{remark}

\section{Two-dimensional complex Jacobian conjecture of degree until 104} \label{degree105-complex}

We consider now the complex case,  {i.e.} considering polynomial maps $F = (F_1, F_2): \C^2 \to \C^2$. 
In this section, we will prove that the complex Jacobian conjecture is true for the degree until 104 
(Theorem \ref{theorem-degree104}). 
For the case of degree 105, we reduce the problem into only four cases 
(Proposition  \ref{prop-degree105}).

Let us recall the principal results concerning the complex Jacobian conjecture under the view of point of Newton polygons studied by Abhyankar \cite{Abhyankar}, Nagata \cite{Nagata2, Nagata3}, Appelgate and Onishi  \cite{Appelgate-Onishi}, Magnus \cite{Magnus}, 
Nakai and Baba \cite{Nakia-Baba} and Moh \cite{Moh}. 
One can also find the results cited below in the book of van den Essen \cite[pages 254-256]{Essen}. 

\begin{theorem} [Abhyankar \cite{Abhyankar}, 1977] \label{thm Abhyankar}
If ${\rm deg} \, (F_1)$ divide  ${\rm deg} \, (F_2) $ or ${\rm deg} \, (F_2) $ divide ${\rm deg} \, (F_1)$,  then the conjecture is true. 
\end{theorem}

Let $d:={\rm gcd} ( {\rm deg}  \, (F_1), {\rm deg} \, (F_2))$ the greatest common divisor of degrees of $F_1$ and $F_2$.

\begin{theorem} [Appelgate and Onishi \cite{Appelgate-Onishi}, Nagata \cite{Nagata2, Nagata3}, 1989]   \label{thm Nataga 1}
  If $d \leq 8$ or $d$ is a prime number then the conjecture is true. 
\end{theorem}

\begin{theorem} [Appelgate, Onishi \cite{Appelgate-Onishi}, Nagata \cite{Nagata3}, 1989] \label{thm Nataga 2} 
If ${\rm deg} \, (F_1)$ or ${\rm deg} \, (F_2)$ is a product of at most two prime numbers, then the conjecture is true. 
\end{theorem}

Recall that ${\rm deg} \, F = {\max} \, \{ {\rm deg} \, (F_1), {\rm deg} \, (F_2)\}$. 
With a highly non-trivial proof, Moh \cite{Moh} proved: 

\begin{theorem} [Moh \cite{Moh}, 1983] \label{thm Moh}
The conjecture is true for ${\rm deg} \, (F) \leq 100$.
\end{theorem}

Using so-called Newton-Puiseux chart, Henryk \.{Z}o\l\c{a}dek \cite{Henryk} give a stronger result than Theorem  \ref{thm Nataga 1}: 

\begin{theorem} [\.{Z}o\l\c{a}dek \cite{Henryk}, 2008] \label{thm Henryk} 
The Jacobian conjecture satisfies for maps with \\
${\rm gcd} ( {\rm deg}  \, (F_1), {\rm deg} \, (F_2)) \leq 16$ and for maps with 
${\rm gcd} ( {\rm deg}  \, (F_1), {\rm deg} \, (F_2))$ equal to 2 times a prime. 
\end{theorem}

In the following, we prove: 

\begin{theorem} \label{theorem-degree104}
The 2-dimensional complex Jacobian conjecture satisfies  for ${\rm deg} \, (F) \leq 104$. 

\end{theorem}

\begin{proof}
At first we remark that if ${\rm deg} \, (F_1) = {\rm deg} \, (F_2)$ then the conjecture is true 
by Theorem \ref{thm Abhyankar}. 
Then we can assume that  ${\rm deg} \, (F_1) < {\rm deg} \, (F_2)$. 

\medskip 

1) The case ${\rm deg} \, (F) =101$: since 101 is a prime, the conjecture is true by Theorem  \ref{thm Nataga 1} or by Theorem \ref{thm Nataga 2}.

\medskip 

2) The case ${\rm deg} \, (F) =102$: 
Assume that 
 ${\rm deg} \, (F_1) < 102$  and ${\rm deg} \, (F_2) = 102$. 
 We see that $102 = 6. 17$ and as $17$ is prime then 
 if ${\rm deg} \, (F_1)$ is not divisible by 17, the greatest common divisor of ${\rm deg} \, (F_1)$ and ${\rm deg } \, (F_2)$ is $d \leq 6$, and the conjecture is true  by Theorem  \ref{thm Nataga 1}. 
 Then we need to consider only the case where ${\rm deg} \, (F_1)$ is divisible by 17, 
 { i.e.} the cases ${\rm deg} (F_1)$ is 17, 34, 51, 68 and 85. 
 If  ${\rm deg} (F_1)$ is 17, 34, 51 and 85, 
 then by Theorem \ref{thm Nataga 2}, 
 the conjecture is true since ${\rm deg} \, (F_1)$ is a product of at most two prime numbers. 
 Now if ${\rm deg} \, (F_1) = 68$, then the greatest common divisor of 
 ${\rm deg} \, (F_1)$ and ${\rm deg} \, (F_2)$ is $34$, which is two times a prime (17 times 2), 
  therefore the conjecture is true by Theorem \ref{thm Henryk}. 

\medskip 

3) The case ${\rm deg} \, (F) =103$: since 103 is a prime, the conjecture is true by Theorem  \ref{thm Nataga 1} or by Theorem \ref{thm Nataga 2}.

\medskip 

4) The case ${\rm deg} \, (F) =104$: 
Assume that ${\rm deg} \, (F_1) < 104$  and ${\rm deg} \, (F_2) = 104$. 
 Since $104 =2^3.13$ then if ${\rm deg} \, (F_1)$ is not divisible by 13, 
 the greatest common divisor of ${\rm deg} \, (F_1)$ and ${\rm deg } \, (F_2)$ is $d \leq 8$, 
 and the conjecture is true  by Theorem  \ref{thm Nataga 1}. 
 Now if ${\rm deg} \, (F_1)$ is not divisible by 2, then  
 the greatest common divisor of ${\rm deg} \, (F_1)$ and ${\rm deg } \, (F_2)$ is $d \leq 13$, 
  the conjecture is also true  by Theorem   \ref{thm Henryk}.  
 Therefore, we need to consider only the case ${\rm deg} \, (F_1)$ is divisible by 26, 
 { i.e.} the cases ${\rm deg} (F_1)$ is $26$, $52$ and $84$. 
 For the two first cases, the conjecture is true by Theorem \ref{thm Abhyankar}. 
 Now if ${\rm deg} \, (F_1) = 84$, the greatest common divisor of 
 ${\rm deg} \, (F_1)$ and ${\rm deg} \, (F_2)$ is $26$, which is two times a prime, 
 then the conjecture is also true by Theorem \ref{thm Henryk}.
\end{proof}

\begin{proposition} \label{prop-degree105}
 If the two-dimensional complex Jacobian conjecture satisfies for the four cases \\ $({\rm deg} \, (F_1), {\rm deg} \, (F_2))$ is  
 (42, 105), (63, 105), (70, 105) and (84, 105) then the conjecture satisfies for polynomial maps of degree 105.
\end{proposition}

\begin{proof}
Assume that 
 ${\rm deg} \, (F_1) < 105$  and ${\rm deg} \, (F_2) = 105$. 
 We see that $105 = 3.5. 7$ and as $7$ is prime then 
 if ${\rm deg} \, (F_1)$ is not divisible by 7, the greatest common divisor of ${\rm deg} \, (F_1)$ and ${\rm deg } \, (F_2)$ is $d \leq 15$, and the conjecture is true  by Theorem \ref{thm Henryk}. 
 Then we need to consider only the case ${\rm deg} \, (F_1)$ is divisible by 7. 
 We remark that by Theorem  \ref{thm Abhyankar} and Theorem \ref{thm Nataga 2}, 
 the cases  we need to consider for ${\rm deg} \, (F_1)$ are only: 7.4, 7.6, 7.8, 7.9, 7.10, 7.12, 7.14. 
 The cases  ${\rm deg} \, (F_1) = 7.4, 7.8$ and 7.14 satisfy 
 the conjecture by Theorem  \ref{thm Nataga 1} since 
 the greatest common divisor of 
 ${\rm deg} \, (F_1)$ and ${\rm deg } \, (F_2)$   in these cases is a prime. 
 So the remain cases are those where ${\rm deg} \, (F_1)$ is 42, 63, 70 and 84.  
 \end{proof}

To end the paper, we would like to mention the two following remarks: 

\begin{remark} \label{Remark-connection}

Theorem \ref{prop1 p:r=1}, Corollary \ref{cor:2 p-r=1} and Theorem \ref{theo n=1} show that 
the complex conjecture is true for the constructed classes 
    $\mathcal{C}^n_{1,s}$ ($n \geq 1$, $s \geq 1$) and 
    $\mathcal{C}^1_{r,s}$ ($r \geq 1$, $s \geq 1$) (see also Remark  \ref{remark classes}). 
    Then for the constructed class $\mathcal{C} = \cup \mathcal{C}^n_{r,s}$, 
    we need to study the conjecture only for the maps belonging to the classes {\rm $\mathcal{C}^n_{r, s}$ ($n \geq 2, r\geq 2, s \geq 1$)} (Corollary   \ref{cor-counterexample}).  
    Using non-proper variables of these classes, one can reduce the degree of the map. 
    For example, if $n = 2$, $r = 2$ and $s= 1$, then the non-proper variables in (\ref{eq:p-r}) are: 
    $$u_0 = y, \quad u_1= x - bx^2y, \quad u_2 = x^2 - b^2x^4y^2.$$
    In this case  ${\rm deg} (u_0) = 1$,  ${\rm deg} (u_1) = 3$ and 
    ${\rm deg} (u_2) = 6$. 
    Then if ${\rm deg}(F) = 105$, each polynomial $F_1$ and $F_2$ is a finite sum of monomials $u_0^{\alpha_0} u_1^{\alpha_1} u_2^{\alpha_2}$ such that 
     $\alpha_0 + 3 \alpha_1 + 6\alpha_2 \leq 105$. \\
   Moreover, with a fixed number $n$ of variables, then by increasing the natural numbers $r$ and $s$,  we can reduce quickly the degree of the map. 
\end{remark}

\begin{remark} 

One can completely generalize the definition of ``non-proper variables" (Definition \ref{d:p-v}) for the general case $F: \K^n \to \K^n$, $n \geq 3$, for both real and complex cases $\K = \R$ and $\K = \C$.
\end{remark} 

\subsection*{\bf Acknowledgement.}  

The author would like to thank Prof. Luis Renato Gon\c calves Dias for many useful discussions  during her visit to the {\it Universidade Federal de Uberl\^andia } (Brazil) in December of 2019. 

This paper was financially supported by the Brazilian Research Support Foundation of Sao Paulo State (FAPESP) process 2023/07802-7 and by the Brazilian National Council for Scientific and Technological Development project CNPq ${\rm n}^{\rm o}$:407454/2023-3. 

\bibliographystyle{plain}

\end{document}